\documentclass [11pt]{amsart}
\usepackage{xcolor,bm}
\usepackage{amssymb,amsmath,amsthm,amsfonts}
\usepackage{enumerate}
\usepackage{setspace}
\usepackage{empheq}
\usepackage[all]{xy}
\usepackage{verbatim}
\usepackage[normalem]{ulem}
\usepackage{todonotes}

\usepackage{hyperref}  
\hypersetup{linkcolor=blue,citecolor=blue,filecolor=blue,urlcolor=blue} 

\numberwithin{equation}{section}
\numberwithin{figure}{section}

\newtheorem{theorem}{Theorem}[section]
\newtheorem{proposition}[theorem]{Proposition}

\theoremstyle{definition}

\usepackage[
backend=biber,
style=alphabetic,
]{biblatex}

\definecolor{myblue}{rgb}{0.6, 0.9, 1}

\makeatletter

\newcommand{\Rmnum}[1]{\expandafter\@slowromancap\romannumeral #1@}
\makeatother

\definecolor{myblue}{rgb}{0.6, 0.9, 1}
\definecolor{mygreen}{rgb}{0,0,1}
\definecolor{purple}{rgb}{0.6,0.2,1}
\definecolor{orange}{rgb}{0.8,0,0.2}

\DeclareMathOperator{\Tr}{Tr}

\newcommand{\eps}{\varepsilon}

\newcommand{\Res}{\operatorname{Res}}

\newcommand{\Gal}{\operatorname{Gal}}

\newcommand{\bb}{\mathbb}

\newcommand{\cal}{\mathcal}
\newcommand{\ovl}{\overline}

\newcommand{\vphi}{\varphi}

\newcommand{\End}{\operatorname{End}}

\newcommand{\Ann}{\operatorname{Ann}}

\newcommand{\disc}{\operatorname{disc}}

\newcommand{\divv}{\operatorname{div}}

\newcommand{\tors}{\operatorname{tors}}
\newcommand{\Frob}{\operatorname{Frob}}
\newcommand{\Art}{\operatorname{Art}}
\newcommand{\cusp}{\operatorname{cusp}}
\newcommand{\ab}{\operatorname{ab}}
\newcommand{\Span}{\operatorname{span}}

\theoremstyle{definition} 

\begin{document}
\title{A New Proof of Merel's Uniform Boundedness Theorem}

\author{Jit Wu Yap}
\address{Department of Mathematics, National University of Singapore, 10 Lower Kent Ridge Road, Singapore 119076}
\email{yapjitwu@nus.edu.sg}

\date{}

\begin{abstract}
We give a new proof of Merel's uniform boundedness theorem for $p$-torsion points of elliptic curves for prime $p$. Our proof constructs holomorphic differentials on $X_0(p)$ from CM points with supersingular reduction at $p$, and show that their $p$-adic Chabuty--Coleman integrals vanish on $J_0(p)(\bb{Q})$. This vanishing follows from studying the action of Frobenius $\Frob_{\ell}$ at primes $\ell$ where an elliptic curve and a fixed CM elliptic curve are geometrically isomorphic after reduction. We then construct a rank zero quotient of the modular Jacobian $J_0(p)$ and verify the formal immersion criterion of Kamienny--Mazur. This new argument was discovered by GPT 6.2 Astra. 
\end{abstract}

	\maketitle

\section{Introduction}
Let $d \geq 1$ be a positive integer and $K$ a number field with $[K:\bb{Q}] \leq d$. Let $E/K$ be an elliptic curve and let $E(K)_{\tors}$ denote the torsion subgroup of $E(K)$, which by Mordell--Weil is a finite set. The following is the famous uniform boundedness theorem of Merel \cite{Mer96}. 

\begin{theorem} \label{IntroTheorem1}
There exists a constant $C = C(d) > 0$, depending only on $d$, such that $|E(K)_{\tors}| \leq C$. 
\end{theorem}

When $d = 1$, this is Mazur's celebrated torsion theorem \cite{Maz77, Maz78}, which even gives a complete list of the possible groups $E(\bb{Q})_{\tors}$ can be. Mazur's strategy is to study rational points on the modular curve $X_0(N)$ and he reduces the problem to showing the existence of a quotient $A$ of the modular Jacobian $J_0(N)$ with $A(\bb{Q})$ having rank zero. Mazur then defines the Eisenstein quotient $\tilde{J}$ and proves that $\tilde{J}(\bb{Q})$ has rank zero, which is the main difficulty in his proof. 
\par 
Kamienny \cite{Kam86} along with Mazur \cite{KM95} extended Mazur's approach to small values of $d$. The case of general $d$ was solved by Merel \cite{Mer96}, by looking at a different quotient of the modular Jacobian $J_0(p)$, called the winding quotient, which we will denote by $J_{e,p}$. The reader may consult \cite{Edi95, Dar04, Reb09} for expositions on the proofs. Similar to Mazur's and Kamienny's approach, Merel requires that $J_{e,p}(\bb{Q})$ has rank zero. This is deduced from deep results of Kolyvagin--Logachev \cite{KL91}, relying on the theory of Euler systems along with the Gross--Zagier formula \cite{GZ86}. This gives an uniform upper bound for torsion of prime order. The remaining case of ruling out large prime power torsion for a fixed prime $p$ can be handled by noting that the gonality of $X_0(p^r)$ tends to $\infty$ \cite{Abr96} and then invoking the Mordell--Lang conjecture which was proven by Faltings \cite{Fal91}, see \cite{Fre94}.
\par 
The aim of this paper is to give a new proof of Theorem \ref{IntroTheorem1} for prime torsion. This new proof was discovered by GPT 6.2 Astra and is relatively elementary, using only basic results coming from $p$-adic modular forms \cite{EMP07}. The general strategy is the same as Kamienny--Mazur \cite{KM95} and Merel \cite{Mer96}, where we want to construct a rank zero quotient of $J_0(p)$ and apply Kamienny's formal immersion criterion. The proof here provides a simpler way to construct such a rank zero quotient. 
\par 
As mentioned in the abstract, the rank zero quotient follows from constructing holomorphic differentials $\omega$ on $X_0(p)$ coming from CM points with supersingular reduction at $p$, whose $p$-adic integral vanish on $J_0(p)(\bb{Q})$. The new fundamental idea is to study primes $\ell$ for our hypothetical elliptic curve $E$ with a rational $p$-torsion point is geometrically isomorphic to a fixed CM elliptic curve after reduction mod $\ell$. 
\par 
Let's explain how one can extract information from such a $\ell$. Such a prime $\ell$ will impose conditions on $E[p^n]$ as a Galois representation, which we can combine with the existence of an order $p$ torsion point. For example, the Weil pairing immediately implies that if each element of $E[p]$ is defined over $\bb{Q}$, then $\zeta_p \in \bb{Q}$ which yields $p \leq 2$. The existence of a rational $p$-torsion point gives us a line in $E[p]$ with every element defined over $\bb{Q}$. The isomorphism $E_{\ell} \simeq E'_{\ell}$ then allows us to read off the action of the Frobenius $\Frob_{\ell}$ on $E[p]$ completely and gain information on $\ell$. We can then vary the CM elliptic curve $E'$ to obtain conditions on $j(E)$. 

\subsection{A simple case} We will first explain a simpler proof that rules out $p$-torsion points over $\bb{Q}$, for all sufficiently large $p \equiv 2 \pmod 3$. More generally, we may replace the condition $p \equiv 2 \pmod 3$ with the set of primes being inert in a fixed imaginary quadratic field $K$. The choice of $K = \bb{Q}(\sqrt{-3})$ gives us $p \equiv 2 \pmod 3$. This proof is quite elementary and does not require any geometry of the modular curve $X_0(p)$, which makes it rather surprising. 
\par 
Let $E$ be our elliptic curve and $P$ be our $p$-torsion point. Let $E'$ be an elliptic curve over $\bb{Q}$ with CM by an order $\mathcal{O}$ in $\bb{Q}(\sqrt{-3})$. Then $j(E') \in \bb{Z}$ and let $\ell$ be a prime dividing $j(E) - j(E')$. If $p$ is sufficiently large, then $E$ must have bad reduction at $3$ and so $\ell \not = 3$. We have the following two cases:
\begin{enumerate}
\item $\ell \equiv 1 \pmod 3$: in this case, $\ell$ splits in $\bb{Q}(\sqrt{-3})$ and hence $E'$ has ordinary reduction mod $\ell$.Then the Frobenius fixes the entire $E[p]$ and the Weil pairing gives us $\ell \equiv 1 \pmod p$.

\item $\ell \equiv 2 \pmod 3$: in this case, $\ell$ is inert in $\bb{Q}(\sqrt{-3})$ and hence $E'$ has supersingular reduction mod $\ell$. Assuming $\ell \geq 5$, we have $E(\bb{F}_{\ell}) = \ell + 1$. Since $P$ specializes to a $p$-torsion point in $E(\bb{F}_{\ell})$, we must have $p \mid \ell + 1$ and so $\ell \equiv - 1 \pmod p$. 
\end{enumerate}

Hence we have 
$$\ell \mid j(E) - j(E') \implies \eps_K(\ell) \ell \equiv 1 \pmod p,$$ where $\eps_K(\ell) = \left(\frac{\ell}{3} \right)$ is the Legendre symbol. We now apply this to two different elliptic curves with CM by $K$. We let the first be given by $y^2 = x^3 + 1$ which has $j$-invariant $0$ and the second by $y^2 = x^3 -15 x + 22$, which has $j$-invariant $54000$. We  define $F(z) = \frac{z - 54000}{z}$. If $j(E)$ was integral at $p$, we would have 
$$F(j(E)) = \frac{j(E) - 54000}{j(E)} \equiv 1 \pmod p.$$
This implies $p \mid 54000$, which can't happen if $p$ is sufficiently large. Hence $p \mid j(E)^{-1}$.  
\par 
The next step is to prove that $p^n \mid j(E)^{-1}$ for all $n \geq 1$, which gives a contradiction. We explain how the proof of $n = 2$ works. The main idea is to apply the comparison with CM elliptic curves on the quotient $E_D = (E/D)$, where $D$ is a cyclic subgroup of $E[p^2]$ containing $P$ and having order $p^2$. The full details are given in Section \ref{sec: p2mod3}. In short, we have a short exact sequence 
$$0 \to V \to E_D[p^2] \to D \to 0.$$
Let $L = \bb{Q}(D)$ be the smallest field for which $E_D$ is defined over. Then $\Gal(\ovl{\bb{Q}}/L)$ acts on both $V$ and $D$ by characters, which we denote by $\psi_V$ and $\psi_D$. Some local calculations and an application of Artin reciprocity to $\frac{\psi_V}{\psi_D}$ and $\left( j(E_D) \right)$ gives us 
$$\prod_{D} F(j(E_D)) \equiv 1 \pmod {p^2}$$
where we take the product over all such $D$'s. We can then show that
$$F(j(E)) \equiv \prod_{D} F(j(E_D)) \equiv 1 \pmod {p^2}$$
giving $p^2 \mid j(E)^{-1}$ as desired. By induction, repeating the same argument gives $p^n \mid j(E)^{-1}$ for all $n \geq 1$, which gives a contradiction. 

\subsection{The general case} In general, the argument in the simpler proof gives us the following: fix an imaginary quadratic field $K$. Then an elliptic curve $E/\bb{Q}$ has no $p$-torsion point for all sufficiently large primes $p$ which are inert $K$. The difficulty in using this to obtain uniformity is that the sufficiently largeness depends on the imaginary quadratic field $K$ itself. 
\par 
To handle the general case, we first geometrize the arguments from before. Let $F(X)$ be a rational function on $X(1) \simeq \bb{P}^1$ with zeroes and poles at CM points which are inert at $p$. Working with the Tate parameter $T$, we have the expansion 
$$j(T) = T^{-1} + 744 + \cdots$$
and we may define $h(T) = F(j(T))$. We let $G_{r}(T) = \prod_{z^{p^r} = T} F(j(z))$ which is a power series in $T$. Letting $r \to \infty$, the power series $G_r(T)$ converges to an actual limit $G_{\infty}(T)$ in the $p$-adic topology. Taking $d \log$, we obtain a formal differential $\omega^{\cusp}_F \in \bb{Z}_p[[T]]dT$. 

\begin{theorem} \label{IntroTheorem2}
There exists a unique $\omega_F \in H^0(X_0(p)_{\bb{Q}_p}, \Omega^1)$ whose restriction to the cusp disk gives us $\omega^{\cusp}_F$. 
\end{theorem}

In fact, $\omega_F$ extends to a section of the relative dualizing sheaf on the integral model $\mathcal{X}_0(p)$. There are a few possible ways to deduce Theorem \ref{IntroTheorem2} and we give one using results of on $p$-adic modular forms \cite{EMP07}. The important property about $F$ is its divisor are CM points with supersingular reduction at $p$. Applying the mod $p^n$ argument in the simpler $p \equiv 2 \pmod 3$ case gives us the following: let $K$ be a number field and let $x = (E,C) \in X_0(p)(K)$. If for all places $v \mid p$ of $K$, the pair $(E,C)$ has potentially multiplicative reduction with $C \simeq \mu_p$ under the Tate isomorphism, then 
\begin{equation} \label{eq: IntroIntegral1}
\sum_{\sigma: K \to \bb{C}_p} \int_{\infty}^{\sigma(x)} \omega_F =0.
\end{equation}
Here, we use the Chabuty--Coleman integral \cite[Section 3]{Sto19}, which defines a pairing on $\Omega^1_J(\bb{Q}_p) \times J(\bb{Q}_p)$. If we let $I_{\omega_F}$ denote the functional $\int_{\infty}^{z} \omega_F$ for $z \in J_0(p)(\bb{Q}_p)$, we can then show that \eqref{eq: IntroIntegral1} implies $I_{\omega_F}$ vanishes on $J_0(p)(\bb{Q})$. 
\par 
This vanishing allows us to construct a quotient of $J_0(p)$ with rank zero over $\bb{Q}$ by taking the ideal $I_F = \{T \in \bb{T} \mid T^* \omega_F = 0\}$ where $\bb{T}$ is the integral Hecke algebra. Using the work of Parent \cite{Par05} which relies on the Zhang's generalization of the Gross--Zagier formula \cite{Zha01}, by varying $F$ we can show that the entire winding quotient $J_{e,p}(\bb{Q})$ has rank zero. This allows us to then directly apply Merel's argument \cite{Mer96} and conclude. 
\par 
To make the argument more self-contained, we will directly show that our quotient satisfies the conditions of the formal-immersion criteria of Kamienny. We first use results of Emerton \cite{Eme02} to reinterpret $I_F$ as Hecke operators that kill an associated divisor $Z_F$ consisting of supersingular elliptic curves mod $p$. We then construct linear functionals on $\bb{T}/(\ell + I_F)\bb{T}$ that kill off $T_1,\ldots,T_m$ but not $T_{m+1}$ mod $\ell$, for $1 \leq m \leq d-1$. This usage of supersingular primes to prove formal imersions has appeared already in the work of Parent \cite{Par05} and Bilu--Parent--Rebolledo \cite{BPR13}.  This handles the case of prime torsion. 
\par 
As remarked earlier, ruling out large prime powers for a fixed prime $p$ is standard and follows easily by combining the Mordell--Lang conjecture, proven by Faltings \cite{Fal91}, with the fact that the gonality of $X_0(p^r)$ goes to $\infty $ as $r \to \infty$ \cite{Abr96}. One interesting thing to note is that Manin \cite{Man69}, following the approach of Demjanenko \cite{Dem66}, already shows that for a fixed number field $K$ and prime $p$, the set $X_0(p^n)(K)$ is finite for sufficiently large $n$. GPT was also able to extend his argument to allow $K$ to vary among number fields with degree $\leq d$, giving a more elementary approach to this case.  

\subsection{Further remarks}
The proofs here were obtained by GPT 6.2 Astra after asking it to use ideas in Alpoge and Claude's proof of the GRH as inspiration for a new proof of Mazur's torsion theorem. The original hope was that since the proof of GRH ports over Stepanov's proof in the case of function fields, perhaps a similar idea would work for other results. In the function field case, the uniform boundedness of torsion for elliptic curves follows easily from $abc$-type statements \cite{HS88, Fre89}, see also \cite{LY26, GG26} for the case of higher dimensional abelian varieties. 
\par 
The proof produced here only took a few prompts and was certainly not what we expected. The inspiration it took from the proof of GRH seems to be the use of CM elliptic curves along with supersingular reduction. GPT 6.2 Astra first produced the proof of $p \equiv 2 \pmod 3$, before finding the proof in the general case. It is interesting to note that before this, we, and also other experts independently, have unsuccessfully tried several times to find a new proof of Merel's theorem using LLMs. It seems that starting with ideas related to the correct one, even if it is related in a vague way, makes these LLMs much more successful at finding the solution. 
\par 
The idea of studying primes $\ell$ dividing $j(E) - j(E')$ for $E'$ a fixed CM elliptic curve has appeared in Elkies' proof \cite{Elk87} of the existence of infinitely many supersingular primes over $\bb{Q}$. It is possible that if asked to use Elkies' proof as inspiration, LLMs would be able to find this proof quickly too. 
\par 
On a side note, effective upper bounds on $C$ were obtained by Merel \cite{Mer96} and Parent \cite{Par99}. Similar effective bounds can be obtained by this proof, although we choose not to work it out in this writeup.

\subsection{AI Disclosure} The original proofs were all obtained by GPT Astra 6.2 with minimal human intervention. The arguments presented here are an outcome after much discussion between the author and GPT in an attempt to digest the proof. Every sentence here was typed by the author and GPT was used to check for mistakes afterwards.

\section{The case of $p \equiv 2 \pmod 3$ over $\bb{Q}$} \label{sec: p2mod3}
We first start off with the case of $p \equiv 2 \pmod 3$ and when $K = \bb{Q}$. This was the first result obtained by GPT 6.2 Astra and is quite elementary. Let $E$ be an elliptic curve over $\bb{Q}$ with a torsion point $P \in E(\bb{Q})$ of order $p$. Let $E'$ be a CM elliptic curve over $\bb{Q}$ with CM by $K = \bb{Q}(\sqrt{-3})$ with order $\mathcal{O}$, so that our prime $p$ is inert in $K$. In particular, $E'$ has supersingular reduction at $p$. We now let $\ell$ be a prime dividing $j(E) - j(E')$. If $p$ is sufficiently large, then $\ell \not = 3$. Note that $\ell \not = p$ too, as otherwise $E$ has supersingular reduction at $p$, implying that $E(\bb{F}_p) = p+1$ which is not compatible with having a rational $p$-torsion point. 

\begin{proposition} \label{Congruence1}
We have the congruence relation 
$$\left(\frac{\ell}{3} \right) \ell \equiv 1 \pmod p$$
where $\left( \frac{\quad \cdot \quad }{\quad \cdot \quad } \right)$ denotes the Legendre symbol.
\end{proposition}

\begin{proof}
 Let's first handle the case where $\ell \equiv 1 \pmod 3$. Then $\ell$ splits in $K$ and so $E'$ has ordinary reduction at $\ell$. Then the action of the Frobenius $\Frob_{\ell}$ can be identified with an element $\pi_{\ell} \in \cal{O}/p\cal{O} \simeq \bb{F}_{p^2}$. However $\pi_{\ell}$ fixes our torsion point $P$, and hence we must have $\pi_{\ell} = 1$ and thus the Frobenius fixes $E[p]$. Since $\det(\Frob_{\ell}) \equiv \ell \pmod p $ \cite[Chapter III, Proposition 8.6]{Sil09}, it follows that $\ell \equiv 1 \pmod p$. 
 \par 
 We now deal with the case where $\ell \equiv 2 \pmod 3$. Then $\ell$ is inert in $K$ and $E'$ has supersingular reduction at $\ell$. If $\ell \geq 5$, by \cite[Chapter V, Exercise 5.10]{Sil09} we have $E(\bb{F}_{\ell}) = \ell + 1$. Since $P$ specializes to a $p$-torsion point of $E(\bb{F}_{\ell})$, we must have $p \mid \ell + 1$ and so $\ell \equiv -1 \pmod p$. This proves our assertion. 
\end{proof}

Now let's consider two different such CM elliptic curves, labelled as $E'_1$ and $E'_2$. We can take two specific examples, one given by $y^2 = x^3 + 1$ with $j$-invariant $0$ and $y^2 = x^3 - 15x + 22$, with $j$-invariant $54000$. Let 
$$F(z) = \frac{z - j(E'_2)}{z - j(E'_1)} = \frac{z - 54000}{z}.$$

\begin{proposition} \label{BadReduction1}
The elliptic curve $E$ has potential multiplicative reduction at $p$, i.e. $p \mid j(E)^{-1}$. 
\end{proposition}

\begin{proof}
Let's assume otherwise that $j(E)$ is integral at $p$. Since
$$F(j(E)) = \frac{j(E)-54000}{j(E)} \equiv 1 \pmod 3,$$
by Proposition \ref{Congruence1}, we also have
$$F(j(E)) = \frac{j(E) - 54000}{j(E)} \equiv 1 \pmod p.$$
As $j(E)$ is integral at $p$, we may multiply over and obtain $54000 \equiv 0 \pmod p$, which is a contradiction for sufficiently large $p$. 
\end{proof}

Proposition \ref{BadReduction1} implies that $p \mid j(E)^{-1}$. We now aim to prove that $p^n \mid j(E)^{-1}$ for all $n \geq 1$. This gives us a contradiction. Since $E$ has potential multiplicative reduction, Tate's uniformization \cite[Chapter V]{Sil94} provides us an isomorphism $E(\bb{C}_p) \simeq \bb{C}_p^{\times}/q^{\bb{Z}}$ for an unique $q \in \bb{Q}_p^{\times}$ with $|q| < 1$. We call $q$ the Tate parameter for $E$. Let $t$ be the Tate parameter for $E/\langle P \rangle$. We will prove the following higher congruence relation. 

\begin{theorem} \label{InductionCongruence1}
For each $r \geq 1$, we have the congruence
$$\prod_{z^{p^r} = t} F(j(z)) \equiv 1 \pmod {p^{r+1}}.$$
Here, we use the $q$-expansion of the $j$-invariant, where 
$$j(z) = z^{-1} + 744 + \cdots$$
so that $j(q) = j(E)$.
\end{theorem}

We first make some observations.

\begin{proposition} \label{CyclicSubgroups1}
Consider the cyclic subgroups $D \subseteq E[p^{r+1}]$ of order $p^{r+1}$ whose order $p$-subgroup is $\langle P \rangle$. Each such $D$ corresponds to a solution of $z^{p^r} = t$, where $E/D$ has Tate parameter $z$. 
\end{proposition}

\begin{proof}
This follows from a easy computation. Indeed, there are $p^r$ different such cyclic subgroups $D$ and each solution of $z^{p^r} = t$ give us such a cyclic subgroup by setting $D = \langle z \rangle$.     
\end{proof}

Now let $E_D$ be the quotient $E/D$ and let $L = \bb{Q}(D)$ be the field of definition for $E_D$. Let $\vphi: E_D \to E$ be the dual isogeny and let $V$ be the kernel. We then have an exact sequence 
$$0 \to V \to E_D[p^{r+1}] \to D \to 0.$$
Since both $V$ and $D$ are cyclic of order $p^{r+1}$, the Galois group $\Gal(\ovl{\bb{Q}}/L)$ acts via characters $\psi_V$ and $\psi_D$, with image landing in $(\bb{Z}/p^{r+1}\bb{Z})^{\times}$. We now apply global Artin reciprocity to the character
$$\theta = \left(\frac{\psi_V}{\psi_D}\right)^{p-1} $$
and the element $a_D = F(j(E_D))$. For each place $v$ of $L$, the local Artin map $\Art_v$ defines a map from $L_v^{\times}$ to $\Gal(L_v^{\ab}/L_v)$ and we have the formula
\begin{equation} \label{eq: Artin1}
\prod_{v} \theta (\Art_v(a_D)) = 1.
\end{equation}

We now compute the local factors at each place $v$. 

\begin{proposition} \label{Artin2}
At a place $v$ above a prime $\ell \not = p$, the local factor $\Art_v(a_D)$ is trivial, i.e. equals $1$. 
\end{proposition}

\begin{proof}
Observe that $\theta$ is valued in a $p$-group, thanks to raising the characters to a $(p-1)$-power. If $j(E_D)$ has a pole at $\ell$, then $F(j(E_D))$ is a principal unit in $O_{L_v}$, which forms a pro-$\ell$ group. Hence $\theta(j(E_D))$ is trivial. 
\par 
Otherwise if $j(E_D)$ is integral, then $E$ has potential good reduction. Since it cannot have additive reduction if $p$ is large enough, it has good reduction and so does $E_D$. Then $\theta$ is unramified and we have the formula
$$\theta(\Art_v(a_D)) = \theta(\Frob_v)^{e_v}$$
where $e_v$ is the exponent of $a_D$ at $v$. If $e_v \not = 0$, then $j(E_D)$ reduces to one of our two CM elliptic curves $E'_1, E'_2$ after reduction mod $v$. We again split into two cases, depending on whether $\ell$ is inert or not in $K$.
\par 
If $\ell$ splits in $K$, then $E'_i$ has ordinary reduction at $\ell$. Since $p$ is inert in $K$, it follows that $E'_i[p^{r+1}]$ is a free $\mathcal{O}/p^{r+1} \mathcal{O}$-module and $\Frob_v$ acts by some element of $\cal{O}/p^{r+1} \cal{O}$. In particular, its action on both $V$ and $D$ are both represented by the same scalar, and hence $\theta(\Frob_v) = 1$. 
\par 
If $\ell$ is inert in $K$, then $E$ has supersingular reduction at $\ell$ and hence the Frobenius at $\ell$ satisfies $\Frob_{\ell}^2 = [-\ell]$. If $f = [k_v:\bb{F}_{\ell}]$ is the inertia degree, then $\Frob_v = \Frob_{\ell}^f$ and so $\Frob_v^2 = [-\ell^f]$. Transporting back to $E_D$, we have that $\Frob_v^2$ acts by a scalar again on $E_D[p^{r+1}]$ and so $\theta(\Frob^2_v) = 1$, which implies $\theta(\Frob_v) = 1$ as it is valued in a $p$-group. 
\end{proof}

We now study what happens for a place $v$ above $p$. 

\begin{proposition} \label{Artin3}
At a place $v$ above the prime $p$, we have the formula
$$\theta(\Art_v(a_D)) \equiv N_{L_v/\bb{Q}_p}(a_D)^{-(p-1)} \pmod {p^{r+1}}.$$
Here, the congruence takes place in $(\bb{Z}/p^{r+1}\bb{Z})^{\times}$.
\end{proposition}

\begin{proof}
Using Tate's uniformization, after possibly passing to a quadratic twist, we see that the kernel $V$ of the dual isogeny $E_D \to E$ is isomorphic to $\mu_{p^{r+1}}$, and that $D$ is isomorphic to the trivial representation $\bb{Z}/p^{r+1} \bb{Z}$. Thus $\theta = \chi_{p^{r+1}}^{p-1}$, where $\chi_{p^{r+1}}$ is the cyclotomic character.
\par 
We now apply norm compatibility \cite[Chapter IV, Proposition 6.2]{Neu99} for local Artin reciprocity, where we use the fields 
$$K = \bb{Q}_p, \quad K' = \bb{Q}(D)_v, \quad L = \bb{Q}_p(\mu_{p^{r+1}}), \quad L' = \bb{Q}(D)_v(\mu_{p^{r+1}})$$ 
in Neukirch's notation. his gives us 
$$\Art_v(a_D) \vert_L = \Art_{\bb{Q}_p}(N_{L_v/\bb{Q}_p}(a_D)) \vert_L.$$
Now the cyclotomic character $\chi_{p^{r+1}}$ sends the RHS to $N_{L_v/\bb{Q}_p}(a_D)^{-1}$ mod $p^{r+1}$ and so we have 
$$\theta(\Art_v(a_D)) \equiv N_{L_v/\bb{Q}_p}(a_D)^{-(p-1)} \pmod {p^{r+1}}$$
as desired.
\end{proof}

Combining Propositions \ref{Artin2} and \ref{Artin3} with \eqref{eq: Artin1}, we obtain 
$$\left(\prod_{v \mid p} N_{L_v/\bb{Q}_p}(a_D) \right)^{-(p-1)} \equiv 1 \pmod {p^{r+1}}.$$
Since each $E_D$ has potential multiplicative reduction, it follows that $a_D$ is a principal unit for $L_v$. Hence we may remove the power and get
\begin{equation} \label{eq: NormCongruence1}
\prod_{v \mid p} N_{L_v/\bb{Q}_p}(a_D) \equiv 1 \pmod {p^{r+1}}.
\end{equation}

\begin{proof}[Proof of Theorem \ref{InductionCongruence1}]
Let $D_1,\ldots,D_{p^r}$ be the cyclic subgroups in Proposition \ref{CyclicSubgroups1}. We split them up into groups of Galois conjugates $S_1,\ldots,S_m$. For each group $S_i$, \eqref{eq: NormCongruence1} implies that 
$$\prod_{D \in S_i} F(j(E_D)) \equiv 1 \pmod {p^{r+1}}.$$
Multiplying them together gives
$$\prod_{i=1}^{p^r} F(j(E_{D_i})) \equiv 1 \pmod {p^{r+1}}.$$
By Proposition \ref{CyclicSubgroups1}, this gives us 
$$\prod_{z^{p^r} = t} F(j(z)) \equiv 1 \pmod {p^{r+1}}$$
as desired. 
\end{proof}

To finish, we note that 
$$\prod_{z^{p^s} = t} F(j(z)) \equiv F(j(t)) \pmod {p^{s+1}}$$
for $s = v_p(t)$. 
Indeed, if we let $h(T) = F(j(T))$, then we have the congruence
\begin{equation} \label{eq: PowerSeriesCongruence1}
\prod_{Z^{p^s} = T} h(Z) \equiv h(T) \pmod { p \bb{Z}_p [[T]]]}.
\end{equation}
Hence if $s = v_p(t)$, Theorem \ref{InductionCongruence1} for $r = s$ with \eqref{eq: PowerSeriesCongruence1} implies that $F(j(t)) \equiv 1 \pmod {p^{s+1}}$. Using the $q$-expansion of $j$, we obtain that $v_p(t) \geq s+1$, giving us a contradiction. Hence for sufficiently large primes $p$ with $p \equiv 2 \pmod 3$, there is no elliptic curve $E$ with a $p$-torsion point $P$ over $\bb{Q}$. 
\par 
In general, the proof works for all sufficiently large primes $p$ that are inert in a fixed imaginary quadratic field $K$. The issue to obtaining uniformity is that the sufficiently large condition depends on the field $K$. 

\section{Constructing holomorphic differentials on $X_0(p)$}
In this section, we will produce holomorphic differentials $\omega$ on $X_0(p)$ whose $p$-adic integral $I_w: J_0(p)(\bb{Q}_p) \to \bb{Q}_p$ vanishes on the rational points $J_0(p)(\bb{Q})$. These differentials will be produced by viewing Theorem \ref{InductionCongruence1} via the geometry on $X_0(p)$. 
\par 
We first start off with a rational function $F$ on $X(1) \simeq \bb{P}^1$ defined over $\bb{Q}$ with zeroes and poles located at CM points with supersingular reduction at $p$. Using the $q$-expansion definition of the $j$-invariant, we define
$$h(T) = F(j(T)) \text{ and } G_r(T) = \prod_{Z^{p^r} = T} h(Z).$$
Note that $h(T) \in 1 +T \bb{Z}_p[[T]]$ since $j(T)^{-1} = T + O(T^2)$ and $F$ has total degree $0$. We may view $G_r(T)$ as the restriction of some global meromorphic function $G_r$ on $X_0(p)$ to the Tate disk at the cusp of $X_0(p)$ where the order $p$ marked subgroup is isomorphic to $\mu_p$. We will denote this cusp by $\infty$.  
\par 
Let $K$ be a number field and let $E$ be an elliptic curve over $K$ and let $D$ be a $\Gal(\ovl{K}/K)$-stable cyclic subgroup such that for every place $v \mid p$, $E$ has potential multiplicative reduction and $D$ is isomorphic to $\mu_p$. Let $t_v$ be the Tate parameter for $E$ at $v$. Then the proof of Theorem \ref{InductionCongruence1} can be extended to show the following theorem. 

\begin{theorem} \label{GeneralCongruence1}
We have 
$$\prod_{v \mid p} N_{K_v/\bb{Q}_p} \left( G_r(t_v) \right) \equiv 1 \pmod {p^{r+1}}.$$
\end{theorem}

Note that $G_r(T)$ is $\Gamma_0(p)$-invariant. Hence if $\omega_{F,r} = d \log G_r(T)$, this gives us formal differentials on the modular curve $X_0(p)$ at the Tate disc at $\infty$. Let $\mathcal{N}$ denote the norm operator $(\mathcal{N} h)(T) = \prod_{Z^p = T} h(Z)$. Then taking $d \log$'s converts $\mathcal{N}$ into the Atkin operator $U_p$ given by
$$U_p \left(\sum_{n \geq 0} a_n T^n \right) = \sum_{n \geq 0} a_{pn} T^n.$$
A simple computation gives us 
$$(d \log)(\mathcal{N} h) = U_p g$$
where $g = d \log h$. Letting $r \to \infty$, we can show that $\omega_{F,r}$ converges $p$-adically to some formal differential $\omega^{\cusp}_F$ that is defined at the Tate disc at $\infty$ for $X_0(p)$. 

\begin{proposition} \label{ConvergenceDifferential}
As $r \to \infty$, the differentials $\omega_{F,r}$ converge $p$-adically uniformly to a formal differential $\omega^{\cusp}_F$. 
\end{proposition}

\begin{proof}
We will prove that $\omega_{F,r+1} \equiv \omega_{F,r} \pmod {p^{r}}$. Since $\omega_{F,r} = d \log G_r$, it suffices to show
$$\log \frac{G_{r+1}(T)}{G_{r}(T)} \equiv 0 \pmod {p^{r}}.$$
We rewrite the LHS as 
$$\log \frac{G_{r+1}(T)}{G_r(T)} = \frac{1}{p} \sum_{S^{p^{r+1}} = T} \log \left(\frac{h(S)^p}{h(S^p)} \right).$$
Now as $h(T) \in 1 + T \bb{Z}_p [[T]]]$, it follows that $h(S)^p \equiv h(S^p) \pmod p$. Hence $h(S)^p/h(S^p)$ expands a series in $S$ with coefficients in $p \bb{Z}_p$ and zero constant term. Summing over $S^{p^{r+1}} = T$ kills off every monomial except for those of the form $S^{p^{r+1}n}$, which contributes $p^{r+1} T^n$. Dividing by $\frac{1}{p}$ gives us the congruence we want. 
\end{proof}

The main result in this section is the following classicality theorem. We remark that GPT also has an alternative proof without using $p$-adic modular forms, which uses \cite[Proposition A.2]{Eme02} instead. 

\begin{theorem} \label{Classicality1}
There exists a global holomorphic differential $\omega_F \in H^0(X_0(p)_{\bb{Q}_p}, \Omega^1)$ whose restriction to the Tate disc at $\infty$ gives us $\omega_F^{\cusp}$. 
\end{theorem}

\begin{proof}
We will use a result of Erickson--Miller--Pixton \cite[Lemma 3.1]{EMP07}, which states that a weight two Serre $p$-adic modular form fixed by $U_p$ is the $q$-expansion of a classical form with $\bb{Q}_p$-coefficients. Clearly $\omega_F$ is fixed by $U_p$. It suffices to check that $\omega_F$ is a weight two Serre $p$-adic modular form. 
\par 
We first check that $d \log h$ gives us a weight two $p$-adic modular form. Since $h = F(j(T))$, it defines a meromorphic weight zero modular form on $X(1)$ with Heegner divisor. Then \cite[Proposition 2.2]{EMP07} gives us that $d \log h$ is a weight two $p$-adic modular form. As $\omega_r$ is obtained by applying $U_p$ repeatedly to $d \log h$, it is also a $p$-adic modular form and hence $\omega_F$ is the $p$-adic uniform limit of $p$-adic modular forms, and hence is a $p$-adic modular form. 
\end{proof}

We can now rephrase Theorem \ref{GeneralCongruence1} in terms of $p$-adic integration over $\omega_F$ via the Chabuty--Coleman pairing \cite[Section 3]{Sto19}. For a curve $C$ over $\bb{Q}_p$ with Jacobian $J$ and origin $O$, The Chabuty--Coleman pairing defines for us a pairing $( \, , \, )$ on $\Omega^1_J(\bb{Q}_p) \times J(\bb{Q}_p) \to \bb{Q}_p$, which we suggestively write as 
$$\int_{O}^{P} \omega := (\omega,P).$$
Now if we fix a basepoint $P_0$ on our curve $C$ so that we have an embedding $C \xhookrightarrow{} J$ and a holomorphic form $\omega$ on $C$, for another point $P$ we may define
$$\int_{P_0}^{P} \omega := \int_{O}^{[P-P_0]} \omega$$
where we view $\omega$ as a $1$-form on $J$. Let $\vphi: D_0 \to C$ be an open disc in $C$, so that $\vphi^* \omega = \omega(z) dz$. We then have the following formula \cite[Section 7]{Sto19}
$$\int_{\vphi(\xi_0)}^{\vphi(\xi_1)} \omega = \int_{\xi_0}^{\xi_1} \omega(z) dz = W(\xi_1) - W(\xi_0)$$
where $W$ is a power series whose derivative is $\omega(z)$. In our situation, if $G_{\infty}$ is the $p$-adic limit of $G_r$, then $d \log G_{\infty} = \omega_{\infty}$ and for $x$ in the Tate disc of the cusp $\infty$ of $X_0(p)$, we have the formula
$$\int_{\infty}^{x} \omega = \log G_{\infty}(x) - \log G_{\infty}(0) = \log G_{\infty}(x).$$
Combining this with Theorem \ref{GeneralCongruence1}, if $P$ is the point defined by $(E,D)$, we obtain 
\begin{equation}  \label{eq: Vanishing1}
    \int_{0}^{R} \omega_F = \sum_{v \mid p} \Tr_{K_v/\bb{Q}_p} \left(\int_{\infty}^{P_v} \omega_F \right) = 0,
    \end{equation}
where $R = \Tr_{K/\bb{Q}}[P-\infty] \in J_0(p)(\bb{Q})$. We now show that as we vary $P$, we obtain all possible $R \in J_0(p)(\bb{Q})$. 

\begin{proposition} \label{Vanishing1}
We have 
$$\int_{0}^{Q} \omega_F = 0$$
for all $Q \in J_0(p)(\bb{Q})$. 
\end{proposition}

\begin{proof}
 For sufficiently large $n$, the Abel-Jacobi map $X_0(p)^{(n)} \to J_0(p)$, where $X_0(p)^{(n)}$ is the symmetric power, given by $Y \mapsto [Y- n (\infty)]$ is a projective bundle by Riemann--Roch and hence smooth. Hence if $D \subseteq X_0(p)$ denotes the Tate disc at $\infty$, the image of $D^n$ in $J_0(p)$ is an open neighbourhood of the origin $O$, say $U$. As $J_0(p)(\bb{Q}_p)$ is compact, there is some $m \geq 1$ such that $mQ \in U$,     
 \par 
 Since $mQ \in U$, there exists some effective divisor $Y$, with no zeroes being $\infty$, defined over $\bb{Q}_p$ for which $[Y - n\infty] = mQ$ in $J_0(p)(\bb{Q}_p)$. Since the space of effective divisors linearly equivalent to $Y$ is a vector space over $\bb{Q}$, we can perturb $Y$ to obtain a linearly equivalent one that is defined over $\bb{Q}$, and such that its zeroes still lie within $D$ and are not $\infty$. 
 \par 
 By \eqref{eq: Vanishing1}, we know that $\int_{O}^{[Y- n (\infty)]} \omega_F = 0$. Since $[Y - n (\infty)] = mQ$, it follows that 
 $$m \int_O^{Q} \omega_F = \int_O^{mQ} \omega_F = 0$$
 as desired. 
\end{proof}

\section{Ruling out large primes}
We are now ready to prove there exists $C = C(d) > 0$, such that for all primes $p > C$ and number fields $K$ with $[K:\bb{Q}] \leq d$, no elliptic curve $E/K$ has a $p$-torsion point in $E(K)$. To do this, we will use Proposition \ref{Vanishing1} to construct a rank zero quotient of $J_0(p)$ and then apply the formal immersion criterion of Kamienny--Mazur. Let's first construct this rank zero quotient. 
\par 
We let $\mathbb{T}$ denote the integral Hecke algebra for $J_0(p)$, i.e. the subring of $\End_{\bb{Q}}(J_0(p))$ that is generated by Hecke operators. For any $T \in \bb{T}$, we have 
\begin{equation} \label{eq: IntegralVanishing1}
\int_{0}^{Q} T^* \omega_F = \int_{0}^{T(Q)} \omega_F = 0
\end{equation}
for all $Q \in J_0(p)(\bb{Q})$. We want to choose the maximal quotient such that Hecke translates of $\omega_F$ span the differential forms in the quotient. Hence we set $I_F = \{T \in \bb{T} \mid T^* \omega_F = 0\}$ and define $A_F = J_0(p)/I_F J_0(p)$. Before proceeding, we will first use a different way to interpret $I_F$, using a result of Emerton \cite{Eme02}. In \cite[Appendix A]{Eme02}, for a suitable unramified extension $K/\bb{Q}_p$, Emerton constructs a residue map
$$\Res: H^0(X_0(p)_{\bb{Q}_p}, \Omega^1) \to M^0 \otimes_{\bb{Z}} K,$$
where $M$ consists of formal sums of supersingular elliptic curves over $\ovl{\bb{F}}_p$ and $M^0$ is the degree zero part.  Emerton \cite[Proposition A.2]{Eme02} shows that this map is Hecke equivariant and is an isomorphism. Hence $I_F$ is also the set of Hecke operators that annihilate $\Res(\omega_F)$. We now calculate $\Res(\omega_F)$. 

\begin{proposition} \label{Calculation1}
Let $Z_F$ be the reduction mod $p$ of $\divv(F)$. Then $\Res(\omega_F) = Z_F$. 
\end{proposition}

\begin{proof}
Let $R$ be the ring of functions on the formal ordinary open $\mathfrak{U}$ of $X_0(p)$ containing the ordinary locus as well as the cusp at $\infty$, i.e. the formal open of points that reduces outside $0$ and the supersingular invariants on the integral model $\mathcal{X}_0(p)$. By restricting to the cusp at $\infty$, we have a map $R \to \bb{Z}_p[[T]]$. This map is injective modulo $p$, i.e. a congruence in $\bb{Z}_p[[T]]$ modulo $p$ implies the same in $R$. 
\par 
The norm operator $\mathcal{N}$ extends to functions of $R$. Hence our power series $G_r(T)$ can be extended to elements of $R$. The proof of Proposition \ref{ConvergenceDifferential} implies that $G_{\infty}(T)$ extends to $R$ too and $G_{\infty}(T)/F(j(T)) \equiv 1 \pmod p$ inside $R$ and thus we can take a $p$-adic logarithm of this function in $R$, say $H$. Hence 
$$\omega_F - d \log F(j(T)) = d \log G_{\infty} - d \log F(j(T)) = dH,$$
which implies that the residues of $\omega_F$ at the supersingular points are the same as that of $F(j(T))$. It is then easy to compute that this gives $\Res(\omega_F) = Z_F$. 
\end{proof}

This gives us an alternative description of $I_F$ by $\{T \in \bb{T} \mid T(Z_F) = 0\}$. The main point here is that $Z_F$ is defined over $\bb{Q}$. An immediate consequence is that
\begin{equation} \label{eq: Annhilator1}
\Ann_{\bb{T} \otimes \bb{Q}_p}(\omega_F) = I_F \otimes_{\bb{Z}} \bb{Q}_p.
\end{equation}

 With \eqref{eq: Annhilator1}, we can now show that $A_F(\bb{Q})$ has rank zero. 

\begin{proposition} \label{RankZero1}
The quotient $A_F(\bb{Q})$ has rank zero.     
\end{proposition}

\begin{proof}
Let $\pi: J_0(p) \to A_F$ be the quotient map. We first claim that 
\begin{equation} \label{eq: DifferentialSpan}
\pi^* H^0(A_{F,\bb{Q}_p}, \Omega^1) = \Span_{\bb{Q}_p} \{T^*\omega_F \mid T \in \bb{T}\}.
\end{equation}
Indeed by construction, $\pi^*H^0(A_{F,\bb{Q}_p}, \Omega^1) = V[I_F]$ where $V = H^0(J_0(p)_{\bb{Q}_p}, \Omega^1)$ and $V[I_F]$ are those that are killed by $I_F$. Clearly $T^*\omega_F$ lie in $V[I_F]$. To prove the opposite inclusion, let $\omega \in V[I_F]$ and let $\lambda$ be a $\bb{Q}_p$-linear functional on $V$ such that $\lambda(T^* \omega_F) = 0$ for all $T \in \bb{T}$. Then there exists some Hecke operator $T \in \bb{T} \otimes_{\bb{Z}} \bb{Q}_p$ such that $\lambda(\omega) = a_1(T \omega)$ where $a_n$ is the $n$th coefficient of its $q$-expansion. Then
$$a_n(T \omega_F) = a_1(T T_n \omega_F) = \lambda(T_n \omega_F) = 0$$
for every $n \geq 1$ and $T \omega_F = 0$. By \eqref{eq: Annhilator1}, it follows that $T \in I_F \otimes_{\bb{Z}} \bb{Q}_p$ and so $\lambda(\omega) = a_1(T \omega) = 0$ as desired, giving us \eqref{eq: DifferentialSpan}
\par 
Using \eqref{eq: IntegralVanishing1}, we obtain that for any $\eta \in H^0(A_{F,\bb{Q}_p}, \Omega^1)$, we have
$$\int_{0}^{\pi(Q)} \eta = 0$$
for all $\pi(Q) \in J_0(p)(\bb{Q})$. This implies that $\pi(Q)$ is torsion \cite[Section 3]{Sto19}. Since $A_F$ is a quotient of $J_0(p)$ over $\bb{Q}$, we have that $J_0(p)(\bb{Q})$ maps to a finite index subgroup of $A_F(\bb{Q})$ and hence $A_F(\bb{Q})$ is torsion and has rank zero as desired.
\end{proof}

We now want to apply the formal immersion criterion of Kamienny \cite{Kam92}. We state it below in the form of \cite[Theorem 5]{Reb09}.

\begin{theorem} \label{Kamienny1} 
Let $p$ be a prime and $\ell \not = p$ a prime. Let $I$ be a saturated ideal of the integral Hecke algebra $\bb{T}$ of $J_0(p)$. Let $A = J_0(p)/I J_0(p)$ and let $\pi: J_0(p) \to A$ be the quotient map. For $d \geq 1$, consider the map 
$$\Phi_d: X_0(p)^d \to A, \quad Y \mapsto \pi([Y - d \infty]).$$
Then the following are equivalent:
\begin{enumerate}
\item $\Phi_{d}$ is a formal immersion when reduced mod $\ell$;
\item There exists $d$ many weight two cusp forms $f_1,\ldots,f_d$ annihilated by $I$ such that the vector $(a_1(f_i),\ldots,a_d(f_i))$ for $i=1,2,\ldots,d$ are linearly independent mod $\ell$.
\item The images of $T_1,\ldots,T_d$ in $\bb{T}/(\ell \bb{T} + I)$ are $\bb{F}_l$-linearly independent.
\end{enumerate}
\end{theorem}

Given a sufficiently large prime $p$, we will choose $F$ and set $I = I_F$. We will then verify that the third condition holds. We will use the description of $I_F = \{T \in \bb{T} \mid T(Z_F) = 0\}$. Recall that $M$ is the space of divisors on supersingular elliptic curves over $\ovl{\bb{F}}_p$. 
\par 
We let $\langle \, , \, \rangle$ denote the standard Euclidean dot product on $M$. Given $x \in M$, we may then define a linear functional on $\bb{T}/(\ell + I_F)\bb{T}$ by $T \mapsto \langle T Z_F, m \rangle \pmod \ell$. We now want to construct $L_x$ so that $L_x(T_i)$ is zero for $1 \leq i \leq k$, but is non-zero for $i = k+1$. Having such $L_x$'s for $k=0,\ldots,d-1$ would immediately imply that $T_1,\ldots,T_d$ in $\bb{T}/(\ell\bb{T}+I)$ are linearly independent. 
\par 
To do so, let $x = [s]$ for some supersingular elliptic curve $s$. we see $L_x(T_k)$ if there are no order $k$ isogeny between $s$ and the support of $Z_F$, and is non-zero if there is, unless there is some cancellation. We thus do the following construction. Let's fix an imaginary quadratic field $K$ for which $p$ is inert, and let $H_1$ be the ring class polynomial of conductor $1$ for $K$, so that $\divv(H_1)$ is the sum of elliptic curves with CM by a maximal order of $K$. Let $A \in M$ be the reduction mod $p$ of $\divv(H_1)$.. We define 
$$F(X) = \left( \displaystyle \prod_{H_1(\alpha) = 0} \Phi_2(X,\alpha) \right)/H_1(X)^3$$
where $\Phi_2$ is the classical modular polynomial. Then the roots of the numerator are exactly those which are related to the roots of $H_1$ by the Hecke correspondence $T_2$. Hence $Z_F = (T_2 - 3)A$.
\par 
We now define $O_{m} = \bb{Z} + m O_K$ and set $E_m = \bb{C}/O_m$. Then $E_m$ has complex multiplication in $K$ and since $p$ is inert, its reduction mod $p$ has supersingular reduction. Let $[s_m] \in M$ represent this reduction and let $L_{s_{2m}} =\langle T Z_F , [s_{2m}] \rangle$.

\begin{proposition} \label{HeckeFunctional1}
Assume that $p > 4Dm^2$ where $D = |\disc(K)|$. We have $L_{s_{2m}}(T_i) = 0$ for all $i < m$ and $L_{s_{2m}}(T_m) \not = 0$ for $\ell$ not dividing a positive integer $b_m \leq \deg(T_m T_2 A)$.  
\end{proposition}

\begin{proof}
We want to read off the coefficient of $s_{2m}$ in $T_i Z_F T = T_i T_2 A - 3 T_i A$. We first apply a result of Kaneko \cite[Theorem 2]{Kan89}. This states that given two CM elliptic curves $E_1, E_2$ with CM by the same field $K$, and whose order has discriminants $D_1,D_2$, then if $D_1 D_2 < p^2$, their reduction mod $p$ must be distinct. This allows us to work with the CM elliptic curves in characteristic zero instead of reducing mod $p$.
\par 
Working over characteristic zero, the lift of $T_i T_2 A$ consists of elliptic curves which are isogenous to one with maximal order by an order $i$ isogeny followed by an order $2$ isogeny and similarly $T_i A$ consists of those for just an order $i$ isogeny. These elliptic curves have conductor dividing $2i < 2m$, and hence their discriminant is at most $(2i)^2 D$. By Kaneko's result, as $(2i)^2 D (2m)^2D < p^2$, their reduction cannot coincide with $s_{2m}$ and hence $L_{s_{2m}}(T_i) = 0$ for all $i < m$. 
\par 
Similarly when $i = m$, we can verify that the coefficient of $s_{2m}$ is zero in $T_i A$. Certainly the coefficient is non-zero in $T_m T_2 A$ as we can reach $s_{2m}$ from a maximal order by a degree $m$ and a degree $2$ isogeny. If we denote the coefficient over $\bb{Z}$ by $b_m$, then $0 < b_m \leq \deg(T_m T_2 A)$ and if $\ell \nmid b_m$, we obtain $L_{s_{2m}}(T_m) \not = 0$.
\end{proof}

To finish, we now apply the theorem of Mazur and Kamienny \cite[Theorem 3]{Reb09}. Combining with Proposition \ref{HeckeFunctional1} and Theorem \ref{Kamienny1}, it suffices to show that for sufficiently large prime $p$, we can find a prime $\ell$ such that $\ell \nmid \prod_{m=1}^{d} b_m$ and such that $p > ( 1 + \ell^{d/2})^2$. We claim that this can be done once $p$ is large enough, depending only on $d$.
\par 
First, by looking at the smallest quadratic nonresidue mod $p$, an elementary argument gives us $K$ such that $p$ is inert and $D_K = O(\sqrt{p})$. In fact with GRH, we can get $O(\log p)$ but $O(\sqrt{p})$ is enough for us. Next, each $b_m$ is bounded above by $\deg(T_m T_2 A) = 3 \sigma_1(m) h_K$, where $h_K$ is the class number and $\sigma_1$ is the sum of divisors. It is easy to prove that $h_K = O(D)$. Hence $\log \left( \prod_{m=1}^{d} b_m \right) = O_d( \log p).$ The product of primes $\ell$ up to some $N$ grows exponentially large and hence once $p$ is sufficiently large, we can always find $\ell$ such that $p > (1 + \ell^{d/2})^2$ and $\ell \nmid \prod_{m=1}^{d} b_m$. By \cite[Theorem 3]{Reb09}, this proves there exists $C = C(d) > 0$, such that there are no elliptic curve $E$ over any number field $L$ with $[L:\bb{Q}] \leq d$ with a $p$-torsion point for $p \geq C$.

\printbibliography

\end{document}